\documentclass[12pt,a4paper]{article}

\usepackage{amsmath,amssymb,amsthm,fancyhdr,supertabular,calc,eepic,fullpage}
\usepackage{mathrsfs}
\usepackage[mathcal]{euscript}
\usepackage{wasysym}

\usepackage[dvips]{graphicx}
\usepackage[dvips]{color}

\newtheorem{theorem}{Theorem}
\newtheorem{lemma}[theorem]{Lemma}

\newtheorem{definition}[theorem]{Definition}
\newtheorem{corollary}[theorem]{Corollary}
\newtheorem*{unsolved}{Unsolved}
\newcommand{\curlyc}[1]{\mathcal{C}_#1}
\newcommand{\trw}{Tr_{GF(q^w)/GF(q)}}
\newcommand{\trtwo}{Tr_{GF(q^2)/GF(q)}}
\newtheorem*{theorema}{Theorem A}
\newtheorem*{theoremb}{Theorem B}
\newtheorem*{theoremc}{Theorem C}

\begin{document}


\pagestyle{plain}

\pagenumbering{arabic}
\title{Polar spaces and embeddings of classical groups}
\author{Nick Gill\thanks{I would like to acknowledge the excellent advice and
support of Associate Professor Tim Penttila.} \\ nickgill@cantab.net}

\maketitle

\begin{abstract}
Given polar spaces $(V,\beta)$ and $(V,Q)$ where $V$ is a vector space over a
field $K$, $\beta$ a reflexive sesquilinear form and $Q$ a quadratic form, we
have associated classical isometry groups. Given a subfield $F$ of $K$ and an
$F$-linear function $L:K\to F$ we can define new spaces $(V,L\beta)$ and
$(V,LQ)$ which are polar spaces over $F$.

The construction so described gives an embedding of the isometry groups of
$(V,\beta)$ and $(V,Q)$ into the isometry groups of $(V,L\beta)$ and $(V,LQ)$.
In the finite field case under certain added restrictions these subgroups are
maximal and form the so called {\it field extension subgroups} of Aschbacher's
class $\curlyc{3}$ \cite{aschbacher}.

We give precise descriptions of the polar spaces so defined and their associated
isometry group embeddings. In the finite field case our results give extra
detail to the account of maximal field extension subgroups given by Kleidman and
Liebeck \cite[p112]{kl}. 

\end{abstract}

\section{Introduction}

Let $(V,\beta)$ and $(V,Q)$ be polar spaces over a field $K$ with $\beta:V\times
V\to K$ a reflexive $\sigma$-sesquilinear form where $\sigma$ is a
$K$-automorphism and $Q:V\to K$ a quadratic form with polar form $f_Q:V\times
V\to K$. Let $F$ be a subfield of $K$ and $L:K\to F$ an $F$-linear function. We
now compose functions to get $L\beta:V\times V\to F$ and $LQ:V\to F$ regarding
$V$ as a vector space over $F$. In order for these forms to be well-defined it
is necessary to impose the condition $\sigma(F)=F$ after which it is easily
verified that $LQ$ is a quadratic form with polar form $Lf_Q$ and $L\beta$ is a
sesquilinear form. In fact if $F\subseteq Fix(\sigma)$ then $\beta$ is
bilinear.

We present three results on this situation: In Section \ref{section:degeneracy},
Theorem A gives conditions on the degeneracy of our composed forms, $L\beta$
and $LQ$. In Section \ref{section:type}, Theorem B gives conditions on the type
(alternating, symmetric or hermitian) of our composed forms. In sections
\ref{section:qtype} and \ref{section:betatype} we consider the situation where our fields are finite.
Theorem C summarises these results and gives the isometry group embeddings which are induced by these
composed forms.

\section{Results on degeneracy}\label{section:degeneracy}

We begin be presenting results on degeneracy. Our definition of degeneracy 
is consistent with that of Taylor \cite{taylor} and so is slightly more 
general than that of Kleidman and Liebeck \cite{kl}:

\begin{definition}
A $\sigma$-sesquilinear form $\beta$ is {\rm non-degenerate} if
$$\beta(u,v)=0, \forall v\in V \implies u=0.$$
A quadratic form $Q$ is {\rm non-degenerate} if its polar form $f_Q$ has the
property that
$$f_Q(u,v)=Q(u)=0, \forall v\in V \implies u=0.$$
The forms are called {\rm degenerate} otherwise.
\end{definition}

Our first result concerns sesquilinear forms and uses an adaptation
of a proof given by Lam \cite{lam}:

\begin{lemma}
$L\beta$ is non-degenerate exactly when $L\neq 0$ and $\beta$ is non-degenerate.
\end{lemma}
\begin{proof}
If either $L=0$ or $\beta$ is degenerate then it is clear that $L\beta$ will be
degenerate. Now suppose that $L\neq 0$, $\beta$ is non-degenerate and $L\beta$
is degenerate. Then there exists nonzero $v\in V$ such that $L\beta(v,w)=0$ for
all $w\in V$. Note that there exists $w\in V$ such that $\beta(v,w)\neq 0$. Now
consider, for any $c\in K$,
\begin{eqnarray*}
\beta(v,\frac{\sigma(c)}{\sigma(\beta(v,w))}w)&=&\sigma(\frac{\sigma(c)}{\sigma(\beta(v,w))})\beta(v,w)
\\
&=&c
\end{eqnarray*}
Then $L\beta(v,\frac{\sigma(c)}{\sigma(\beta(v,w))}w)=Lc=o$ for all $c\in K$.
This implies that $L=0$ which is a contradiction.
\end{proof}

We turn our attention to quadratic forms. To begin with we can apply the
previous lemma directly to get the following:

\begin{lemma}
If $L=0$ or $Q$ is degenerate then $LQ$ is degenerate. If $f_Q$ is
non-degenerate then $LQ$ is non-degenerate.
\end{lemma}

Thus we are left with the question of what happens when $Q$ is non-degenerate
and $f_Q$ is degenerate. This can only occur in characteristic 2. We are able
to present results only for the case where $V$ is finite-dimensional and $K$ is finite, in which
case we have the following well-known result (see, for example \cite[p.
143]{taylor}):

\begin{theorem}
A non-degenerate quadratic form $Q$ on a vector space $V$ over $GF(2^h)$ has a
degenerate associated polar form if and only if dim $V$ is odd, in
which case the radical of $f_Q$, $rad(V,f_Q)$, is of dimension 1.
\end{theorem}

\begin{corollary}
Let $K$ be a finite field of characteristic 2. Suppose $dim_KV$ is odd, $f_Q$
is degenerate and $Q$ is non-degenerate. Then $LQ$ is degenerate.
\end{corollary}
\begin{proof}
Take $x\in rad(V,f_Q)$. Then $x\in rad(V,Lf_Q)$. Hence $rad(V,Lf_Q)\supseteq
rad(V,f_Q)$. But $dim_F(rad(V,Lf_Q))\geq dim_F(rad(V,f_Q))>1.$ Hence $LQ$ is
degenerate.
\end{proof}

We can summarise our main results in the following:

\begin{theorema}
Let $\beta:V\times V\to K$ be a reflexive $\sigma$-sesqulinear form. Let
$Q:V\to K$ be a quadratic form. let $F$ be a subfield of $K$ and $L:K\to F$ be
a $F$-linear function. Then:
\begin{itemize}
\item $L\beta$ is non-degenerate if and only if $\beta$ is non-degenerate and
$L\neq 0$;
\item If $char \ K \neq 2$, or $K=GF(2^h)$ for some integer $h$ and $dim_KV$ is even, then $LQ$ is
non-degenerate if and only if $Q$ is non-degenerate and $L\neq 0$;
\item If $K=GF(2^h)$ for some integer $h$ and $dim_KV$ is odd then $LQ$ is degenerate;
\end{itemize}
\end{theorema}

\begin{unsolved}
We have failed to ascertain the conditions under which $LQ$ is degenerate
in the case where $char \ K=2$, $|K|+dim_KV$ is infinite, $Q$ is non-degenerate and
$f_Q$ is degenerate.
\end{unsolved}

\section{A classification of $\beta$ into form}\label{section:type}

Taking reflexive sesquilinear form $\beta:V\times V\to K$ to be alternating,
symmetric or hermition, $L:K\to F$, $F$-linear and not identically zero, we
seek to classify $L\beta$ into these three categories or else as being
`atypical', i.e. not of of this form. 

The conditions under which $\beta$ is hermitian, $char \ K = 2$ and $L\beta$ is
alternating will prove to be the most difficult and we discuss this case first. 
Observe that we must have $F\subseteq Fix(\sigma)$.

Let $\sigma$ be the field automorphism of order $2$ associated with $\beta$. It
is easily shown that $K/Fix(\sigma)$ is a Galois extension and we may
therefore define a trace function:
$$Tr_{K/Fix(\sigma)}:K\to Fix(\sigma), x\mapsto x+\sigma(x).$$
Now any $Fix(\sigma)$-function $L:K\to Fix(\sigma)$ can be written in the
form, for some $\alpha\in K$,
$$L:K\to Fix(\sigma), x\mapsto Tr_{K/Fix(\sigma)}(\alpha x).$$

\begin{lemma}
When $char \ K = 2$ and $\beta$ is hermitian, $L\beta$ is alternating if and
only if $F\subseteq Fix(\sigma)$ and $L\sigma = L$.
\end{lemma}
\begin{proof}
Write $L:K\to F,x\mapsto L_1\circ Tr_{K/Fix(\sigma)}(\alpha x)$ for some
$\alpha\in K$ and some $L_1:Fix(\sigma)\to F$, $F$-linear and not identically
zero. We suppose that 
$Tr_{K/Fix(\sigma)}(\alpha\sigma)=Tr_{K/Fix(\sigma)}(\alpha)$ and it is enough to prove that 
$Tr_{K/Fix(\sigma)}(\alpha\beta)$ is alternating. Now for $x\in K,$
\begin{eqnarray*}
Tr_{K/Fix(\sigma)}(\alpha\sigma(x))=Tr_{K/Fix(\sigma)}(\alpha x) &\implies &
\alpha\sigma(x)+\sigma(\alpha\sigma(x)) = \alpha x+ \sigma (\alpha x) \\
&\implies& \alpha\sigma(x) +\sigma(\alpha)x=\alpha x + \sigma(\alpha)\sigma(x)
\\
&\implies& (\sigma(\alpha)+\alpha)(\sigma(x)+x)=0.
\end{eqnarray*}

Since $\sigma(x)+x\neq 0$ for all $x\not\in Fix(\sigma)$, we must have
$\sigma(\alpha)=\alpha$. Then
\begin{eqnarray*}
Tr_{K/Fix(\sigma)}(\alpha\beta)(x,x)&=&\alpha\beta(x,x)+\sigma(\alpha\beta(x,x)) \\
&=& \alpha\beta(x,x)+ \sigma(\alpha)\sigma\beta(x,x) \\
&=& (\alpha+\sigma(\alpha))\beta(x,x)=0.
\end{eqnarray*}

\end{proof}

We are now able to state our main result:

\begin{theoremb}
Let $\beta:V\times V\to K$ be a reflexive sesquilinear form. Let $K/F$ be a
field extension with $L:K\to F$ a $F$-linear function which is not identically
zero. Then we classify $\beta$ into type as follows:
\begin{itemize}

\item	If $\beta$ is alternating then $L\beta$ is alternating;

\item	If $\beta$ is symmetric then $L\beta$ is symmetric;

\item	If $char \ K=2$, $K$ is finite and $\beta$ is symmetric not alternating then $L\beta$ is
symmetric not alternating;

\item If $\beta$ is hermitian  and $F\not\subseteq Fix(\sigma)$ then
\begin{enumerate}
\item $L\beta$ is hermitian if and only if $L\sigma = \sigma L$;
\item $L\beta$ is atypical if and only if $L\sigma \neq \sigma L$;
\end{enumerate}
\item If $\beta$ is hermitian and $F\subseteq Fix(\sigma)$ then
\begin{enumerate}
\item $L\beta$ is symmetric if and only if $L\sigma = L$;
\item $L\beta$ is alternating if and only if $char \ K \neq 2$ and $L\sigma=-L$
OR $char \ K = 2$ and $L\sigma = L$;
\item $L\beta$ is atypical if and only if $L\sigma\neq\pm L$.
\end{enumerate}

\end{itemize}
\end{theoremb}
\begin{proof}
The first two statements are self-evident. 

We turn to the third statement. Given $\beta$ symmetric not alternating,
$L\beta$ will be alternating if and only if $\{\beta(x,x)| x\in V\}\subseteq
null(L)$. Since $L\neq 0$ it is enough to show that $f:V\to K, x\to
\beta(x,x)$ is onto. Take any $x\in V$ such that $\beta(x,x)=a\in K^*$. Take
any $c\in K$. Then $\beta(\sqrt{\frac{c}{a}}x,\sqrt{\frac{c}{a}}x)=c$ as
required.

For the remainder we assume that $\beta$ is hermitian. First of all suppose
that $F\not\subseteq Fix(\sigma)$ so $L\beta$ is $\sigma$-sesquilinear. Then
$L\beta(v_1,v_2) = L\sigma\beta(v_2,v_1)$ for any $v_1,v_2\in V$ and so $L\beta$
is hermitian if and only if $L\sigma|_{\Im(\beta)}=\sigma L|_{\Im(\beta)}$.
Since $\beta$ is surjective we are done.

Next suppose that $F\subseteq Fix(\sigma)$ in which case 
$L\beta(v_1,v_2)=L\sigma\beta(v_2,v_1)$. This is symmetric if and only if
$L\sigma|_{\Im(\beta)}=L|_{\Im(\beta)}$ and so $L\beta$ is symmetric exactly when
$L\sigma = L$.

Now we examine when $L\beta$ is alternating. When $char \ K$ is odd this is
equivalent to $L\beta$ being skew-symmetric which, by an analagous argument to
the symmetric case, occurs exactly when $L\sigma = -L$. When $char \ K = 2$ the
previous lemma gives us the required result. The only other possibility is for
$L\beta$ to be atypical hence we have our final equivalence.
\end{proof}

\begin{unsolved}
We have failed to ascertain the conditions under which $L\beta$ is alternating
in the case where $char \ K=2$, $K$ is infinite and $\beta$ is symmetric not
alternating.
\end{unsolved}

\section{The isometry classes of $(V,LQ)$ over finite
fields}\label{section:qtype}

Define $Q:V\to GF(q^w)$ a non-degenerate quadratic form, $L:GF(q^w)\to GF(q)$ a $GF(q)$-linear function which is not the zero function and $Tr_{GF(q^w)/GF(q)}$ the trace function. We restrict $V$ to be a finite $A$-dimensional vector space over $GF(q^w)$. In order to classify $(V,Q)$ into isometry classes we need to examine the situation when $Aw$ is even and distinguish between the $O^+$ and $O^-$ cases.

Our first lemma will be useful in distinguishing the isometry class of $LQ$ as
well as giving an application of the classification:

\begin{lemma}
The isometry group for $Q$, $Isom(Q,V)$, is a subgroup of the isometry group
for $LQ$, $Isom(LQ,V)$.
\end{lemma}
\begin{proof}
Simply observe that if $T:V\to V$ satisfies $Q(Tu)=Q(u)$ for all $u\in V$ then
$LQ(Tu)=LQ(u)$ for all $u$ in $V$.
\end{proof}

Consider first the situation when $A$ is even:

\begin{lemma}\label{lemma:oplus}
Let $(V,q)$ have isometry class $O^+(A,q^w)$. Then $(V,LQ)$ has isometry class
$O^+(Aw,q)$. Thus $O^+(A,q^w)\leq O^+(Aw,q)$.
\end{lemma}
\begin{proof}
Let $W$ be a maximal totally singular subspace of $(V,q)$. Then
$dim_KW=\frac{1}{2}dim_KV$. But $W$ is also a totally singular subspace of
$(V,LQ)$ and $dim_FW=\frac{1}{2}dim_FV$. Thus $(V,LQ)$ is of type $O^+(Aw,q)$.
\end{proof}

\begin{lemma}\label{lemma:ominus}
Let $(V,q)$ have isometry class $O^-(A,q^w)$. Then $(V,LQ)$ has isometry class
$O^-(Aw,q)$. Thus $O^-(A,q^w)\leq O^-(Aw,q)$.
\end{lemma}
\begin{proof}
{\bf Suppose first of all that $A=2$}. Suppose in addition that $(V,LQ)$ has
isometry class $O^+(2w,q)$. Then $O^+(2,q^w)\leq O^-(Aw,q)$ and so, by the
theorem of Lagrange,
$$(q^w+1)\big| q^{w(w-1)}\prod^{w-1}_{i=1}(q^{2i}-1).$$
If a primitive prime divisor of $q^{2w}-1$ exists then this is impossible hence
we must deal with the exceptions given by Zsigmondy. The first possibility is
that $w=1$, in which case $L:GF(q^w)\to GF(q^w)$ has form $x\mapsto ax$ for
some $a\in GF(q^w)^*$. Clearly an element of $V$ is singular under $Q$ exactly
when it is singular under $LQ$. Then $(V,q)$ and $(V,LQ)$ have the same Witt
index and hence share type which is a contradiction.

The second possibility is that $(q,w)=(2,3)$ in which case we must consider
whether or not $O^-(2,8)\leq O^+(6,2)$. Examining the atlas \cite{atlas} we see that
$O^-(2,8)$ contains elements of order $9$ while $O^+(6,2)$ does not, hence this
possibility can be excluded.

{\bf Now suppose that $A= 2m+2$ for some $m\geq 1$.} Then $V=U\perp W$ under $Q$
where $U$ is a direct sum of $m$ hyperbolic lines and $W$ is an anisotropic
subspace of dimension 2. Then $Q|_U$ is of type $O^+$ and hence $LQ|_U$ is of
type $O^+$. Similarly $Q|_W$ is of type $O^-$ and, since $dim_{GF(q^w)}W=2$,
$LQ|_W$ is of type $O^-$. Then $V=U\perp W$ under $LQ$, $U$ is a direct sum of
$mw$ hyperbolic lines under $LQ$ and $W$ is a direct sum of $w-1$ hyperbolic
lines with a 2-dimensional anisotropic subspace under $LQ$. Hence $(V,LQ)$ is
of type $O^-$. 
\end{proof}

We now consider the situation when $A$ is odd. If the characteristic equals 2 then Theorem A implies that $LQ$ is degenerate so we exclude this situation. We will be interested in the situation where $w$ is even and the characteristic is odd. We will write $L:GF(q^w)\to GF(q)$ in the form, for some $\alpha\in GF(q^w)^*$, 
$$L=Tr_{GF(q^w)/GF(q)}(\alpha):GF(q^w)\to GF(q), x\mapsto \sum^{w-1}_{i=0}(\alpha x)^{q^i}.$$
We will need to work with the {\it discriminant} of our form $LQ$ for which we
will need two preliminary results:

\begin{lemma}
Let $q$ be odd, $k\in GF(q^2)\backslash GF(q)$ such that $k^2\in GF(q)$. Then
$$Tr_{GF(q^2)/GF(q)}(k)=0.$$
\end{lemma}
\begin{proof}
Observe that $GF(q^2)=GF(q)(k)$ and $k$ has minimum polynomial $x^2-k^2$. Now $Gal(GF(q^2)/GF(q))$ acts on the set of roots of this minimum polynomial. Since the trace map is the sum of the elements of $Gal(GF(q^2)/GF(q))$, $Tr_{GF(q^2)/GF(q)}(k)=k-k=0$.
\end{proof}

\begin{theorem}\label{theorem:discriminant}
A non-degenerate quadratic form, $Q:V\to GF(q)$ where $V$ is a $2n$-dimensional
vector space over $GF(q)$ and $q$ is odd, gives rise to an $O^+(2n,q)$ space if
and only if $disc(Q)\equiv(-1)^n(mod \ GF(q)^{*2})$. Here $GF(q)^{*2}$ is the subgroup of
$GF(q)^*$ consisting of all square terms.
\end{theorem}

A proof of the previous theorem can be found, for instance, in
\cite[p.32]{kl}. We can now proceed with our study of the type of $LQ$.

\begin{lemma}
Let $(V,Q)$ be of type $O$ and be one-dimensional over $GF(q^2)$, $q$ odd. Then $Q$
has form $Q(u)=\gamma u^2$, for some $\gamma\in GF(q^w)^*$. Then $LQ$ has type,
\begin{eqnarray*}
O^+	&\iff&(\alpha\gamma)^{-2}\in GF(q)\backslash GF(q)^{*2} \ or \
(\alpha\gamma)^{q+1}\not\equiv-1(mod \ GF(q)^{*2}), \\ 
O^- 	&\iff&(\alpha\gamma)^{-2}\not\in GF(q)\backslash GF(q)^{*2} \ and \
(\alpha\gamma)^{q+1}\equiv-1(mod \ GF(q)^{*2}).
\end{eqnarray*}
\end{lemma}
\begin{proof}
Observe that $LQ:V\to GF(q), u\mapsto \alpha\gamma u^2 + (\alpha\gamma u^2)^q$ has polar
form $Lf_Q:V\times V\to GF(q), (u,v)\mapsto u^TMv$ where, over a basis for $V$ over
$GF(q)$, $\{1,\omega\},$
$$M=\left(\begin{array}{cc}
2\trtwo(\alpha\gamma) & 2\trtwo(\alpha\gamma\omega) \\
2\trtwo(\alpha\gamma\omega) & 2\trtwo(\alpha\gamma\omega^2)
\end{array}
\right).$$

Now take $f$ to be an element of $GF(q)$ such that $\sqrt{f}\not\in GF(q)$. Then 
$(\alpha\gamma)^{-2}\in GF(q)\backslash GF(q)^{*2}$ if and only if $(\alpha\gamma)^{-1}\sqrt{f}\in GF(q)$. 

{\bf Suppose that $(\alpha\gamma)^{-1}\sqrt{f}\not\in GF(q)$}. Let
$\omega=(\alpha\gamma)^{-1}\sqrt{f}$. Then
$$M=\left(\begin{array}{cc}
2\trtwo(\alpha\gamma) & 2\trtwo(\sqrt{f}) \\
2\trtwo(\sqrt{f}) & 2\trtwo(\alpha^{-1}\gamma^{-1}f)
\end{array}
\right).$$
Then the discriminant of the form $LQ$ is
\begin{eqnarray*}
&&4\trtwo(\alpha\gamma)\trtwo(\alpha^{-1}\gamma^{-1}f)-4(\trtwo(\sqrt{f}))^2 \\
&=& 4f(\alpha\gamma+(\alpha\gamma)^q)((\alpha\gamma)^{-1}+(\alpha\gamma)^{-q})\
\ \ \ \ {\rm (since} \ \trtwo(\sqrt{f})=0) \\
&=&\frac{f2^2(\trtwo(\alpha\gamma))^2}{(\alpha\gamma)^{q+1}}.
\end{eqnarray*}
Referring to Theorem \ref{theorem:discriminant} we see that our result holds in
this case.

{\bf Suppose that $(\alpha\gamma)^{-1}\sqrt{f}\in GF(q)$}. Let
$\omega=\sqrt{f}$. Then $\alpha\gamma=f_2\sqrt{f}$ for some $f_2\in GF(q)$ and
$$M=\left(\begin{array}{cc}
2\trtwo(f_2\sqrt{f}) & 2\trtwo(f_2f) \\
2\trtwo(f_2f) & 2\trtwo(f_2f\sqrt{f})
\end{array}
\right).$$
The discriminant of the form $LQ$ is
\begin{eqnarray*}
&&4\trtwo(f_2\sqrt{f})\trtwo(f_2f\sqrt{f})-4(\trtwo(f_2f))^2 \\
&=& -4(\trtwo(f_2f))^2 \ \ \ \ \ {\rm (since} \ \trtwo(\sqrt{f})=0). 
\end{eqnarray*}
Appealing to Theorem \ref{theorem:discriminant} we conclude that $LQ$ is of
isometry class $O^+$ in all cases here.



\end{proof}

\begin{lemma}
Let $(V,Q)$ be $A$-dimensional of type $O$ over field $GF(q^w)$ of odd
characteristic. Let $S$ be a non-dimensional anisotropic subspace (or {\it
germ}) where $Q\big|_S$ has form $Q(s)=\gamma s^2$ for some $\gamma\in
GF(q^w)^*$. Let $w=2n$. Then $LQ$ has type 
\begin{eqnarray*}
O^+	&\iff&(\alpha\gamma)^{-2}\in GF(q^n)\backslash GF(q^n)^{*2} \ or \
(\alpha\gamma)^{q+1}\not\equiv-1(mod \ GF(q^n)^{*2}),  \\
O^- 	&\iff&(\alpha\gamma)^{-2}\not\in GF(q^n)\backslash GF(q^n)^{*2} \ and \
(\alpha\gamma)^{q+1}\equiv-1(mod \ GF(q^n)^{*2}).
\end{eqnarray*}
\end{lemma}
\begin{proof}
First take $A$ odd and $w=2$. Then $(V,Q)=(R,Q\big|_S)\perp(S,Q\big|_R)$ where $R$ is an
orthogonal direct sum of orthogonal hyperbolic lines and $S$ is a
one-dimensional anisotropic orthogonal space. Then $LQ\big|_R$ is of type $O^+$
and $LQ\big|_S$ will be either of type $O^+$ or $O^-$ according to the
conditions of the previous lemma. Since
$(V,LQ)=(R,LQ\big|_R)\perp(S,LQ\big|_S)$ the type of $LQ$ is determined
according to the conditions given.

Now take $A$ odd, $w$ any even number. Then
$L=Tr_{GF(q^n)/GF(q)}\circ Tr_{GF(q^w)/GF(q^n)}\circ K$ where $K:GF(q^w)\to
GF(q^w), x\mapsto \alpha x$. By the previous paragraph the conditions of the
theorem are the conditions under which $Tr_{GF(q^w)/GF(q^n)}\circ K\circ Q$ will
be of type $O^+$ or $O^-$. By Lemmas \ref{lemma:ominus} and \ref{lemma:oplus}
we know that further compositions with $ Tr_{GF(q^n)/GF(q)}$ will not change
this type. The result follows.

\end{proof}

We will summarise the results of this section and the next in Theorem C at the
end of the paper.

\section{The isometry classes of $(V,L\beta)$ over finite
fields}\label{section:betatype}

Define $\beta:V\times V\to GF(q^w)$ to be a non-degenerate reflexive sesquilinear form of one of the three types, $V$ $A$-dimensional over $GF(q^w)$. Define $L:GF(q^w)\to GF(q)$, $GF(q)$-linear and not the zero function.

If we consider $\beta$ symmetric over a field of odd characteristic then
$\beta$ shares isometry class with the quadratic form
$Q(v)=\frac{1}{2}\beta(v,v)$ and the results of the previous section give the type of
$L\beta$. 

Similarly if $\beta$ is alternating or if the characteristic is 2 and $\beta$ is symmetric
not alternating, then Theorem B gives the type of $L\beta$. Note that over
finite fields, symmetric not alternating forms result in polar spaces which are
called {\it pseudo-symplectic}.

In this section we need to consider the case where $\beta$ is hermitian with
automorphism $\sigma$. Once
again we will take $L=\trw(\alpha)$ for some $\alpha$ in $GF(q)^*$. Consider
first the case where $GF(q)\not\subseteq Fix(\sigma)$ which occurs exactly when
$w$ is odd:

\begin{lemma}
Let $\beta$ be hermitian. When $w$ is odd,
\begin{enumerate}
\item	$L\beta$ is hermitian $\iff$ $\sigma(\alpha)=\alpha$;
\item	$L\beta$ is atypical $\iff$ $\sigma(\alpha)\neq\alpha$.
\end{enumerate}
When $w$ is even,
\begin{enumerate}
\item	$L\beta$ is symmetric $\iff$ $\sigma(\alpha)=\alpha$;
\item	$L\beta$ is alternating $\iff$ $\sigma(\alpha)=-\alpha$;
\item	$L\beta$ is atypical $\iff$ $\sigma(\alpha)\neq\pm\alpha$.
\end{enumerate}
\end{lemma}
\begin{proof}
Observe that $L\beta$ is bilinear if and only if $F(q)\subseteq Fix(\sigma)$ if and only if $w$ is even.

Suppose first that $w$ is odd; then it is enough to prove the first
equivalence. By Theorem B we know that $L\beta$ is hermitian if and only if
$L\sigma = \sigma L$. Now
$$L\sigma(x) = \sigma L(x) \iff \sigma\trw((\sigma(\alpha)-\alpha)x)=0.$$
The surjectivity of the trace function gives us our result.

Now suppose that $w$ is even. By Theorem B it is enough to prove that
$L\sigma=\pm L \iff \sigma(\alpha)=\pm\alpha$. Let the Galois group of the field extension
$GF(q^w)/GF(q)=\{\sigma_1,\dots,\sigma_w\}$. Then
\begin{eqnarray*}
L\sigma=\pm L&\iff& \sum_{i=1}^w\sigma_i(\alpha\sigma(x)) =
\pm \sum_{i=1}^w\sigma_i(\alpha x) \\
&\iff& \sum_{i=1}^w\sigma_i\sigma(\alpha\sigma(x)) =
\pm \sum_{i=1}^w\sigma_i(\alpha x) \\
&\iff& \sum_{i=1}^w\sigma_i((\sigma(\alpha)\mp\alpha)x) = 0.
\end{eqnarray*}
Once again the surjectivity of the trace function gives us our result.
\end{proof}

To complete the classification we need to ascertain the isometry group of
$L\beta$ in the case where it is symmetric.

\begin{lemma}
Suppose that $\beta$ is hermitian, $w$ is even, $\sigma(\alpha)=\alpha$ and $V$
is $A$-dimensional over $GF(q^w)$. Then the isometry class of $L\beta$ is,\
$$O^+(Aw,q)\iff A \ is \ even, \ O^-(Aw,q)\iff A \ is \ odd.$$
\end{lemma}
\begin{proof}
If $A$ is even then, with respect to the hermitian form $\beta$, $V$ contains a totally isotropic
subspace of dimension $\frac{A}{2}.$ This subspace is also totally isotropic
with respect to $L\beta$ and over $GF(q)$ has dimension $\frac{Aw}{2}$. Hence
$(V,L\beta)$ is of type $O^+$.

Now take $A$ to be odd. First suppose that $A=1$ and $w=2$ so that
$GF(q)=Fix(\sigma)$. Then, given a basis for $V$ over $GF(q)$, $\{1,\omega\}$,
we have $\beta(x,y)=x\sigma(y)$, $L(x)=\alpha x+\sigma(\alpha x)$ and the
matrix of $L\beta$ is
$$\left(\begin{array}{cc}
2\trw(\alpha) & 2\trw(\alpha\omega) \\
2\trw(\alpha\omega) & 2\trw(\alpha\omega\sigma(\omega))
\end{array}
\right).$$
Now put $\omega=\sqrt{f}$ where $GF(q^2)=GF(q)(\sqrt{f})$ and the discriminant
of $L\beta$ is $-4f\alpha^2$. Since this is minus a non-square, $L\beta$ is of
type $O^-$.

Now let $A$ be any odd integer, $w=2$. Then
$(V,\beta)=(R,\beta\big|_S)\perp(S,\beta\big|_R)$ where $R$ is an
orthogonal direct sum of orthogonal hyperbolic lines and $S$ is a
one-dimensional unitary space. Then $L\beta\big|_R$ is of type $O^+$ by the
first part of this lemma, $L\beta\big|_S$ is of type $O^-$ by the previous
argument and hence $(V,L\beta)$ is of type $O^-$.

Finally suppose that $w>2$, in which case $L =
Tr_{GF(q^{\frac{w}{2}})/GF(q)}\circ Tr_{GF(q^w)/GF(q^{\frac{w}{2}})}(\alpha)$.
We know that $Tr_{GF(q^w)/GF(q^{\frac{w}{2}})}(\alpha\beta)$ is of type $O^-$;
then Lemma \ref{lemma:ominus} implies that $(V,L\beta)$ is of type $O^-$.
\end{proof}

We are now in a position to summarise the results of the last two sections.

\begin{theoremc}
Let $V$ be an $A$-dimensional polar space over $GF(q^w)$. Take $L: GF(q^w)\to
GF(q),x\mapsto \trw(\alpha x)$ for some $\alpha\in GF(q^w)^*$.

Suppose first of all that $V$ is defined via a quadratic form $Q:V\to GF(q^w)$.
If the form has a germ $U$ then  $Q\big|_U(x)=\gamma x^2$ for some $\gamma\in GF(q^w)^*$. Then we classify $LQ$ into type, including the classical group embedding,as follows:
\begin{center}
\begin{tabular}{|c|c|c|c|}
\hline 
\hline 
Type of $Q$ & Type of $LQ$ & Conditions & Embedding \\
\hline 
\hline 
$O^+$ & $O^+$ & always & $O^+(A,q^w)\leq O^+(Aw,q)$ \\
\hline 
$O^-$ & $O^-$ & always & $O^-(A,q^w)\leq O^-(Aw,q)$ \\
\hline 
O & degenerate & $q$ even & - \\
\hline 
$O$ & $O$ & $w$ odd, $q$ odd & $O(A,q^w)\leq O(Aw,q)$ \\
\hline 
$O$ & $O^+$ & $w$ even, $q$ odd; & $O(A,q^w)\leq O^+(Aw,q)$ \\
&& $(\alpha\gamma)^{-2}\in GF(q^{\frac{w}{2}})\backslash GF(q^{\frac{w}{2}})^{*2}  \ or$  & \\
&& $(\alpha\gamma)^{q+1}\not\equiv-1(mod \ GF(q^{\frac{w}{2}})^{*2})$ & \\
\hline 
$O$ & $O^-$ & $w$ even, $q$ odd; & $O(A,q^w)\leq O^-(Aw,q)$ \\
&& $(\alpha\gamma)^{-2}\not\in GF(q^{\frac{w}{2}})\backslash GF(q^{\frac{w}{2}})^{*2}  \ and$ & \\ 
&& $(\alpha\gamma)^{q+1}\equiv-1(mod \ GF(q^{\frac{w}{2}})^{*2})$ & \\
\hline 
\hline 
\end{tabular}
\end{center}

Suppose next that $V$ is defined via a reflexive $\sigma$-sesquilinear
form $\beta:V\times V\to GF(q^w).$  If the
characteristic is odd and $\beta$ is symmetric then the type of $L\beta$ and
its associated classical group embedding is given in the previous table taking
$Q$ to be the quadratic form $Q(v)=\frac{1}{2}\beta(v,v)$.

In all other cases the type of $L\beta$, with associated classical group
embedding, is as follows:

\begin{center}
\begin{tabular}[c]{|c|c|c|c|}
\hline 
\hline 
Type of $\beta$ & Type of $L\beta$ & Conditions & Embedding \\
\hline 
\hline 
hermitian & hermitian & $w$ odd, $\sigma(\alpha)=\alpha$ & $U(A,q^w)\leq
U(Aw,q)$ \\
\hline 
hermitian & atypical & $w$ odd, $\sigma(\alpha)\neq\alpha$ & - \\
\hline 
hermitian & alternating & $w$ even, $q$ even, $\sigma(\alpha)=\alpha$ &
$U(A,q^w)\leq
Sp(Aw,q)$ \\
\hline 
hermitian & alternating & $w$ even, $q$ odd, $\sigma(\alpha)=-\alpha$ &
$U(A,q^w)\leq
Sp(Aw,q)$ \\
\hline 
hermitian & atypical & $w$ even, $\sigma(\alpha)\neq\pm\alpha$ & - \\
\hline 
hermitian & $O^+$ & $w$ even, $q$ odd, $A$ even, $\sigma(\alpha)=\alpha$ &
$U(A,q^w)\leq
O^+(Aw,q)$ \\
\hline 
hermitian & $O^-$ & $w$ even, $q$ odd, $A$ odd, $\sigma(\alpha)=\alpha$ &
$U(A,q^w)\leq
O^-(Aw,q)$ \\
\hline 
alternating & alternating & always & $Sp(A,q^w)\leq
Sp(Aw,q)$ \\
\hline 
pseudo & pseudo & $q$ even & - \\
-symplectic & -symplectic & & \\
\hline
\hline 
\end{tabular}
\end{center}
\end{theoremc}

\bibliographystyle{plain}
\bibliography{classicalpaper}

\end{document}